\documentclass[a4paper, 10pt, twocolumn]{article}

\pdfoutput=1

\usepackage{amsfonts,amssymb,amsmath,amsthm}
\usepackage{subfigure}
\usepackage{graphicx}
\usepackage[footnotesize]{caption}
\usepackage{mystyle}
\usepackage{url}
\graphicspath{{./images/}}

\usepackage[top=1.5cm, left=1.5cm, right=1.5cm, bottom=1.5cm]{geometry}

\def\R{\mathbb{R}}

\def\N{\mathbb{N}}

\def\closure{\mathrm{closure}}
\def\ext{\mathrm{ext}}
\def\rext{\mathrm{rext}}

\def\M{\mathcal{M}}

\newtheorem{assumption}{Assumption}

\newtheorem{theorem}{Theorem}
\title{Convex Regularization and Representer Theorems}

\author{C. Boyer$^1$, A. Chambolle$^2$, Y. De Castro$^{3,4}$, V. Duval$^{4,5}$, F. de Gournay$^6$ and P. Weiss$^7$.\\
\footnotesize $^1$ LPSM, Sorbonne Universit\'e, ENS Paris
\footnotesize $^2$ CMAP, CNRS, Polytechnique, Palaiseau \\
\footnotesize $^3$ LMO, Universit\'e Paris-Sud, Orsay 
\footnotesize $^4$ MOKAPLAN, INRIA, Paris \\
\footnotesize $^5$ Ceremade, Universit\'e Paris Dauphine
\footnotesize $^6$ INSA, Universit\'e de Toulouse
\footnotesize $^7$ ITAV, CNRS, Universit\'e de Toulouse
}
\date{\empty} % no need for a date

\renewenvironment{abstract}{\bf\small {\em\ Abstract---}}{}

\begin{document}

\maketitle

\begin{abstract} 
We establish a result which states that regularizing an inverse problem with the gauge of a convex set $C$ yields solutions which are linear combinations of a few extreme points or elements of the extreme rays of $C$. These can be understood as the \textit{atoms} of the regularizer. We then explicit that general principle by using a few popular applications. In particular, we relate it to the common wisdom that total gradient variation minimization favors the reconstruction of piecewise constant images. 
 \end{abstract}

\newcommand{\myparagraph}[1]{\vspace{2mm}\noindent\textbf{#1}}

%%%%%%%%%%%%%%%%%%%%%%%%%%%%%%%%%%%%%%%%%%%%%%%%%%%%%%%%%%%%%%%%%%
%%%%%%%%%%%%%%%%%%%%%%%%%%%%%%%%%%%%%%%%%%%%%%%%%%%%%%%%%%%%%%%%%%
%%%%%%%%%%%%%%%%%%%%%%%%%%%%%%%%%%%%%%%%%%%%%%%%%%%%%%%%%%%%%%%%%%
\section{Introduction}

Let $E$ denote a Hausdorff locally convex vector space, and $m\in \N$. Let $\Phi : E\to \R^m$ be a bounded linear mapping called \emph{sensing operator} and $u\in E$ denote a signal. 
We assume that a device returns a set of measurements of the form $y=P(\Phi u)$, where $P:\R^m\to \R^m$ is a perturbation operator, such as quantization (1-bit compressed sensing), modulus (phase retrieval), additive Gaussian noise,... 
We consider the inverse problem of estimating $u$ from $y$. 

An approach at the heart of many successful approach consists of solving problems of the form:
\begin{equation}\label{eq::mainproblem}
 \inf_{u \in E} f(\Phi u) + J_C(u), 
\end{equation}
where $J_{C}$ is a \emph{gauge function} associated to a convex set $C \subset E$ defined by 
\begin{equation}
 J_{C}(u)=\inf \{ \lambda \text{ s.t. } u\in \lambda C, \lambda \ge 0\},
\end{equation}
and $f$ is an arbitrary convex or non-convex function called \emph{data fitting term}. This function should depend on the choice of the data $y$ and the perturbation $P$. Interpreting $J_{C}$ as the \emph{atomic norm} corresponding to a set of so-called \emph{atoms}, it suggested in~\cite{chandrasekaran2012convex} that $J_{C}$ favors solutions $u$ which consist of a few such atoms.

The present note makes that idea precise by describing the faces (in particular the extreme points) of the solution set of~\eqref{eq::mainproblem}. 
When the minimizer is unique, our theorem describes the solution as a conical combination of less than $m$ extreme points or elements of the extreme rays of $C$.

We then showcase the strength of this theorem by applying it to popular regularizers such as total variation (mass), nonnegativity constraints, positive semi-definite constraints or total gradient variation. 
This note is a reduced version of a longer work by the same authors, which contains all the proofs and additional results \cite{BCCDGW}.

\section{Definitions and preliminaries}

\paragraph{Lines, half-lines, linear closeness.}
A \textit{line} is an affine subspace of $E$ with dimension $1$. 
An open half-line is a set $\rho$ of the form $\rho=\{p+tv, t>0\}$, where $p,v\in E$, $v\neq 0$. 
An open half-line $\rho$ contained in $C$ is called a \emph{ray} of $C$. 

\paragraph{Lineality space.}
The \emph{lineality} space $C_K$ of a closed convex set $C\subseteq E$ is a subspace defined by $C_K\eqdef\{u\in E, u+C=C\}$.
Accordingly, $C$ can be decomposed as:
\begin{equation}\label{eq:decomposition}
C=C_K + C_B, \mbox{ where } \mbox{ and } C_B=C\cap W, 
\end{equation}
where $W$ is a linear complement of $C_K$ in $E$. The set $C_B$ contains no line.
There is a parallel to make between the set $C_K$ and the kernel of a linear operator, justifying the notation $C_K$.

\paragraph{Extreme points and extreme rays.}
An \textit{extreme point} of $C$ is a point $p\in C$ such that $C\setminus \{p\}$ is convex. 
We let $\ext(C)$ denote the set of extreme points of $C$.
An \textit{extremal ray} of $C$ is a ray $\rho \in C$ such that if $x,y\in C$ and the open segment $(x,y)$ intersects $\rho$, then $(x,y)\subset \rho$. The set of extreme rays of $C$ is denoted $\rext(C)$.

The following result generalizes the Krein-Milman theorem:
\begin{theorem}[Klee 1957 \cite{klee_extremal_1957} \label{thm:klee}]
 If $C\in E$ is locally compact, closed, convex and contains no line, then
 \begin{equation}
  C = \closure(\mathrm{conv}(\ext(C) \cup \rext(C))).
 \end{equation}
\end{theorem}

\paragraph{Faces.}
Following~\cite{klee_theorem_1963}, if $p\in C$, the smallest face of $C$  which contains $p$ is denoted $\face{p}{C}$. It is defined as the union of $\{p\}$ and all the open segments in $C$ which have $p$ as an inner point. The dimension of $\face{p}{C}$ is defined as the dimension of its affine hull. The collection of all elementary faces, $\{\face{p}{C}\}_{p\in C}$, is a partition of $C$.
Extreme points correspond to the zero-dimensional faces of $C$, while extreme rays are subset of the one-dimensional faces. 

\section{A representer theorem}

In what follows, we let $t^\star=J_C(u^\star)$ with $u^\star \in S^\star$ be the infimum in \eqref{eq::mainproblem} and we define
\begin{equation}
\delta= \begin{cases}1 &\textrm{ if } t^\star =\inf_{u\in E} J_C(u), \\ 0&\textrm{ if } t^\star >\inf_{u\in E} J_C(u).\end{cases}
\end{equation}

We decompose $S^\star$ as $S^\star = S_K^\star + S_B^\star$ following \eqref{eq:decomposition} and let $d$ be the dimension of $\Phi(C_K)$.

Our main results will hold under the following assumptions. 
\begin{assumption}[Main assumptions]\label{assumptions} \ 
\begin{itemize}
 \item $C$ is nonempty, closed and convex,
 \item $C_B$ is closed, locally compact.
 \item $S^\star$ is not empty. 
\end{itemize}
\end{assumption}

Let us state a result valid for arbitrary data fitting functions $f$.
\begin{theorem}
\label{thm1}
Assume that Assumption \ref{assumptions} is satisfied. 
Then there exists at least one solution $u^\star\in S^\star$ of the form:
\begin{equation}
 u^\star = u_K^\star + \sum_{k=1}^r \alpha_k \psi_k,
\end{equation}
where
\begin{itemize}
 \item $u_K^\star \in C_K$, and $\alpha_k\geq 0$ for all $k$.
 \item $r\leq m-d+\delta$ and $\psi_k \in \ext(C_B)$ for all $k$,
 \item or $r\leq m-1-d+\delta$ and $\psi_k \in \rext(C_B)\cup \ext(C_B)$ for all $k$.
\end{itemize}
\end{theorem}

When $f$ is a convex, lower semi-continuous function, a stronger result can be obtained.
\begin{theorem}
\label{thm2}
Assume that Assumption \ref{assumptions} is satisfied and that $f$ is convex and closed.%\footnote{Ce r\'esultat n'apparait dans le papier (je ne fait pas d'hypoth\`ese de forte convexit\'e sur $f$), mais je crois bien qu'il est vrai. L'id\'ee est que si $J_C$ varie sur $S^\star$ (ce qui est d\'ej\`a peu probable), ca ne change pas les points extr\'emaux. C'est un des int\'er\^ets des jauges outre la lecture plus facile.} 

Then $S_K^\star=C_K\cap \ker(\Phi)$. Moreover, for every $p\in S_B^\star$, let $j$ be the dimension of the face $\face{p}{S_B^\star}$. Then $p$ can be written as a conical combination of:
 \begin{itemize}
  \item $m+j-d+\delta$ extreme points of $C_B$,
  \item or $m+j-1-d+\delta$ points of $C_B$, each an extreme point of $C$ or in an extreme ray of $C$.
 \end{itemize}
\end{theorem}
This result actually describes the structure of the \emph{whole} solution set. Indeed, by taking $j=0$ and $j=1$, we describe the extreme points and the extreme rays of the solution set, which are enough to reconstruct $S^\star$ entirely thanks to Theorem \ref{thm:klee}.

\section{Examples of applications}

\subsection{$\ell^1$-norm and total variation}

The $\ell^1$-norm $J_C(u) = \|u\|_1$ on $E=\R^n$ is the gauge of the unit $\ell^1$-ball $C=\{u\in \R^n, \|u\|_1\leq 1\}$. 
It is well known that $\ext(C) = \{\pm e_i, 1\leq i \leq m\}$, where $e_i$ denote the $i$-th element of the canonical basis. In addition, $\rext(C)=\emptyset$. Hence, our theorems state that some solutions of \eqref{eq::mainproblem} will be $m$-sparse, whatever the data fitting term $f$. This result is one of the (implicit) cornerstones in compressed sensing \cite{candes2006robust}.

Similarly, let $E=\M$ denote the space of Radon measures on a domain $\Omega$ of $\R^d$.
The total variation in the measure-theoretic sense (or mass) of a measure $u\in \M$ can be written as $J_C(u)=\|u\|_\M$, where $C=\{u \in \M, \|u\|_\M\leq 1\}$. We have $\ext(C)=\{\pm \delta_x, x\in \Omega\}$ and $\rext(C)=\emptyset$. This explains why using the total-variation as a regularizer allows to recover sparse spikes \cite{candes2012super,duval2013exact}.

\subsection{Generalized total variation}

Let $D'(\Omega)$ denote the space of distributions, and $E\subset D'(\Omega)$. 
Let $L:E\to \M$ denote a surjective continuous linear operator, $C\eqdef\{u\in E, \|Lu\|_\M\leq 1\}$ and the associated regularizer $J_C(u) = \|Lu\|_\M$. Then the Fisher-Jerome theorem \cite{fisher1975spline}, which was recently revisited in \cite{unser2017splines,flinth2017exact} can be obtained using our theorems and by remarking that $C_K=\ker(L)$, that $\ext(C_B)=\{\pm L^+\delta_x, x\in \Omega\}$ and that $\rext(C)=\emptyset$. This theorem essentially states that the solutions of \eqref{eq::mainproblem} can be expressed as a small linear combination of splines.

\subsection{Nonnegativity}

Set $E=\R^n$ and define $J_C(u)=\chi_C(u)$, where $C$ is the nonnegative orthant. Then $\ext(E)=\{0\}$ and $\rext(C) = \{\alpha e_i, \alpha\geq 0, 1\leq i \leq n\}$. Hence our main theorems state that at least one solution of \eqref{eq::mainproblem} is $m$-sparse. This result helps understanding the field of nonnegative least-squares for instance. It is also a critical (and implicit) element of \cite{donoho2005sparse}.

Similarly, let $E=\R^{m\times m}$ and define $C$ as the cone of positive symmetric semi-definite matrices. Our main theorem in that case states that some solutions will be rank $m$ matrices, which is a key fact to  understand the success of matrix completion \cite{candes-tao1}.

\subsection{Total gradient variation}

The representation principle can be applied with the total variation in the functional analysis sense, i.e. the total variation of the gradient \cite{rudin1992nonlinear}. Given the set $E\eqdef L^{d/(d-1)}(\R^d)$, we define $J_C$ as the (isotropic) total gradient variation, \ie{} 
\begin{align*}
& J_C(u) \eqdef \|Du\|_{\M^d} \\
&= \sup\left(\int u\mathrm{div}(\phi) \, dx, \phi \in C^1_c(\R^d)^d, \sup_{x\in \R^d} \|\phi(x)\|_2\leq 1\right),
\end{align*}
where $D$ is the distributional gradient and $\|\cdot\|_{\M^d}$ denotes a vectorial total variation. If $F\subset\RR^d$ has finite measure, we define its perimeter as $P(F)\eqdef J_C(\bun_{F})$.

The following result is due to Fleming~\cite{fleming1957functions}, with a more exhaustive treatment of the problem covered in~\cite{ambrosio2001connected}. We refer to those references for the proper definition of simple sets.

\begin{prop}[\cite{fleming1957functions,ambrosio2001connected}]
The extreme points of  
\begin{align}
  C= \{u\in L^{d/(d-1)}(\RR^d), \|Du\|_{\M^d}\leq 1 \},
\end{align}
 are the functions $u=\pm\bun_{F}/P(F)$, where $F$ is a simple set and $P(F)<+\infty$. 
\end{prop}
Informally, simple sets are the simply connected sets in the measure-theoretic sense (\ie{} they consist of one connected component and they have no holes). 

Combining this result with our theorems allows us to conclude that at least one solution of total gradient variation regularized problems are the sum of $m$ indicator functions. 
This result is a clear explanation of the staircasing effect when using a finite number of measurements $m$.
In addition, it gives some insight on the family of functions that can be exactly recovered by total gradient variation minimization. %Surprisingly, no satisfactory exact reconstruction guarantees are available in the literature, despite the good empirical performance of this prior. 
We believe that this result is an important step towards a better understanding of this regularizer.

\newpage
%%%%%%%%%%%%%%%%%%%%%%%%%%%%%%%%%%%%%%%%%%%%%%%%%%%%%%%%%%%%%%%%%%
%%%%%%%%%%%%%%%%%%%%%%%%%%%%%%%%%%%%%%%%%%%%%%%%%%%%%%%%%%%%%%%%%%
%%%%%%%%%%%%%%%%%%%%%%%%%%%%%%%%%%%%%%%%%%%%%%%%%%%%%%%%%%%%%%%%%%
%% You can make the bibliography smaller
\bibliographystyle{abbrv}
{\small
\bibliography{biblio}
}

\end{document}